\documentclass[a4paper,10pt]{report}

\usepackage{graphicx}
\usepackage{amsfonts}
\usepackage{titlesec}
\usepackage{amsmath,amssymb,amsthm}
\usepackage{times}
\usepackage[left=2.5cm,top=3cm,bottom=3cm,right=2.5cm]{geometry}

\usepackage{fancyhdr}
\titlespacing\subsection{0pt}{\parskip}{-\parskip}

\newcommand{\R}{{\mathbb R}}

\newcommand{\N}{{\mathbb N}}

\newtheorem{re}{Remark}

\newtheorem{cor}{Corollary}

\newtheorem{lem}{Lemma}

\newtheorem{theo}{Theorem}

\newcommand{\ds}{\displaystyle}

\newcommand{\la}{\lambda}

\newcommand{\ga}{\gamma}
\newcommand{\al}{\alpha}

\newcommand{\su}{\subset}

\newcommand{\lt}{\left\|}
\newcommand{\rt}{\right\|}
\newcommand{\ltl}{\left\langle}
\newcommand{\rtl}{\right\rangle}

\newcommand{\ra}{\rightarrow}
\newcommand{\lra}{\longrightarrow}

\newenvironment{prf}{\emph{\textbf{Proof. }}}{\hfill\rule{0.5em}{0.5em}}

\date{\textbf{SEPTEMBER 2014}}

\begin{document}
\thispagestyle{plain}

\pagenumbering{arabic}
\setcounter{page}{1}

\begin{center} \Large HYBRID ALGORITHM FOR NONLINEAR EQUILIBRIUM, \\ VARIATIONAL INEQUALITY AND FIXED POINT PROBLEMS.
\end{center}

 \vspace{0.5cm}

 \begin{center}
\normalsize O. I. Agha Ibiam,$^{1}$ \hspace{0.2cm} L. O. Madu, $^{2}$ \hspace{0.2cm} E. U. Ofoedu,$^{3,*}$ \hspace{0.2cm} C. E. Onyi $^{4}$ \hspace{0.2cm} and \hspace{0.2cm} H. Zegeye$^{5}$
\end{center}

 \begin{center}
\footnotesize $^{1,}$ $^{2}$ $^{3}$ $^{4}$ Department of Mathematics, Nnamdi Azikiwe University, P.M.B. 5025, Awka, Anambra State, Nigeria.\\
$^5$Department of Mathematics, University of Botswana, Gaborone, Botswana.
 \end{center}

 \vspace{1cm}

\subsection*{\normalsize ABSTRACT} 
\noindent In this paper, we introduce an iterative process which converges strongly to a common element of sets of solutions of finite family of generalized equilibrium problems, sets of fixed points of finite family of continuous relatively nonexpansive mappings and sets of zeros of finite family of $\ga$-inverse strongly monotone mappings in Banach spaces. Our theorems improve and generalize several results which are announced recently. Our iteration process, method of proof and corollaries are of independent interest.

\setcounter{section}{1}
\vspace{0.5cm}
\noindent
2010 {\itshape Mathematics Subject Classification} 47H05, 47H06, 47H10, 47J05, 47J25.

\vspace{0.25cm}
\noindent
\textbf{KEYWORDS}: Equilibrium problems, generalized projection, relatively quasi-nonexpansive mappings, monotone mappings, variational inequality problems.

\subsection*{\normalsize INTRODUCTION}
\normalsize
\noindent
Let $E$ be a real Banach space with dual $E^*$. We denote by $J$ the normalized duality mapping from $E$ to $2^{E^*}$ defined by $$Jx := \{f^* \in E^* : \ltl x, f^*\rtl = \lt x\rt^2 = \lt f^*\rt^2\},$$ where $\ltl \cdot, \cdot\rtl$ denotes the generalized duality pairing between members of $E$ and members of $E^*$. It is well known that if $E^*$ is strictly convex then $J$ is single valued and if $E$ is uniformly smooth then $J$ is uniformly continuous on bounded subsets of $E$. Moreover, if $E$ is reflexive and strictly convex Banach space with a strictly convex dual, then $J^{-1}$ is single valued, one-to-one, surjective, and it is the duality mapping from $E^*$ into $E$ and thus $JJ^{-1} = I_{E^*}$ and $J^{-1}J = I_E$ (see e.g., \cite{cioranescu, takahashi1}).

\begin{re}
If $E = l_p, \; 1<p<\infty,$ then $J$ is given explicitly as $Jx = \lt x\rt^{2-p}_{l_p}y \in l_q, \; x = \{x_1, x_2, \cdots\}$, $y = \{x_1|x_1|^{p-2}, x_2|x_2|^{p-2}, \cdots\}$, $\frac{1}{p}+\frac{1}{q} = 1$ and so $J^{-1}$ is also easily computed. Furthermore, if $E = L_p, \; Jx = \lt x\rt^{2-p}_{L_p}|x|^{p-2}x \in L_q$, $\frac{1}{p}+\frac{1}{q} = 1$ (see e.g., \cite{alber}). Observe that if $H$ is a Hilbert space, then $E^* = H$, $J$ and $J^{-1}$ becomes the identity operator on $H$.

\noindent Let $C$ be a nonempty closed convex subset of a real Banach space $E$. Let $f : C \times C \ra \R$ be a bifunction, where $\R$ is the set of real numbers, and $B: C \ra E^*$ be a nonlinear mapping. The generalized equilibrium problem (for short, $GEP$) for $f$ and $B$ is to find $u \in C$ such that
\begin{equation}
f(u, v)+\ltl Bu, v-u\rtl \geq 0, \; \forall \; v \in C.
\end{equation}
\end{re}

\let\footnotesize\scriptsize
\let\thefootnote\relax\footnote{\textbf{$^*$Corresponding Author}:  E-mail addresses: eu.ofoedu@unizik.edu.ng, euofoedu@yahoo.com, euofoedu@gmail.com}

\fancyhead[l]{\scriptsize \textbf{\itshape Agha Ibiam, \hspace{0.1cm} Madu, \hspace{0.1cm} Ofoedu, \hspace{0.1cm} Onyi \hspace{0.1cm} and \hspace{0.1cm} Zegeye}}
\fancyhead[r]{\scriptsize \textbf{\itshape Journal of Physical Research, ISSN: 2141-8403 prints, Vol 7, No. 1, Jan 2017}}
\renewcommand{\headrulewidth}{0.1mm}

\noindent The set of solutions for the problem (1) is denoted by $GEP(f, B)$, i.e., $$GEP(f, B) := \{u \in C: f(u, v)+\ltl Bu, v-u\rtl \geq 0, \; \forall \; v \in C\}.$$ If $B=0$ in (1), then $GEP(1)$ reduces to the classical equilibrium problem (for short, $EP$) and $GEP(f, 0)$ is denoted by $EP(f)$, i.e., $$EP(f) = \{u \in C: f(u, v) \geq 0, \; \forall \; v \in C\}.$$
If $f=0$ in (1), then $GEP(1)$ reduces to the classical variational inequality problem and $GEP(0, B)$ is denoted by $VI(B, C)$, i.e., $$VI(B, C) = \{u^* \in C : \ltl Bu^*, v-u^*\rtl \geq 0, \; \forall \; v \in C\}.$$

\noindent The problem (1) is very general in the sense that it includes, as special cases, optimization problems, variational inequalities, min-max problems, the Nash equilibrium problems in non-cooperative games and countless of others; (see, for instance, \cite{blum, moudafi})

\noindent A mapping $A: D(A) \su E \ra E^*$, is said to be monotone if for each $x, y \in D(A)$, the following inequality holds:
\begin{equation}
\ltl x-y, Ax-Ay\rtl \geq 0.
\end{equation}
A is said to be $\ga$-\emph{inverse strongly monotone} if there exists a positive real number $\ga$ such that
\begin{equation}
\ltl x-y, Ax-Ay\rtl \geq \ga\lt Ax-Ay\rt^2, \quad \mbox{for all} \quad x, y \in K.
\end{equation}
If $A$ is $\ga$-inverse strongly monotone, then it is \emph{Lipschitz continuous} with constant $\frac{1}{\ga}$, that is, $\lt Ax-Ay\rt \leq \frac{1}{\ga}\lt x-y\rt,$ for all $x, y \in D(A)$, and hence uniformly continuous.

\noindent Let $E$ be a smooth real Banach space. The function $\phi: E \times E \ra \R$ defined by
\begin{equation}
\phi(x, y) = \lt x\rt^2-2\ltl x, Jy\rtl+\lt y\rt^2 \quad \mbox{for all} \quad x, y \in E,
\end{equation}
is studied by Alber \cite{alber}, Kamimura and Takahashi \cite{kamimura}, and Riech \cite{Reich}. It is obvious from the definition of the function $\phi$ that
\begin{equation}
(\lt x\rt - \lt y\rt)^2 \leq \phi(x, y) \leq (\lt x\rt+\lt y\rt)^2 \quad \mbox{for} \quad x, y \in E.
\end{equation}
Observe that in a Hilbert space $H$, (5) reduces to $\phi(x, y) = \lt x-y\rt^2,$ for $x, y \in H.$

\noindent Let $E$ be a reflexive, strictly convex and smooth real Banach space and let $C$ be a nonempty closed and convex subset of $E$. The \emph{generalized projection mapping}, introduced by Alber \cite{alber}, is a mapping $\Pi_C : E \ra C$, that assigns to an arbitrary point $y \in E$ the minimum point of the functional $\phi(\cdot, y),$ that is $\Pi_C y = \bar{x}$, where $\bar{x}$ is the solution to the minimization problem 
\begin{equation}
\phi(\bar{x}, y) = \min\{\phi(x, y), x \in C\}.
\end{equation}
In fact, we have the following result.
\begin{lem}\cite{alber}
Let $C$ be a nonempty closed and convex subset of a reflexive, strictly convex, and smooth real Banach space $E$ and let $y \in E.$ Then there exists a unique element $x_0 \in C$ such that $\phi(x_0, y) = \min\{\phi(z, y): z \in C\}.$
\end{lem}

\noindent
Let $C$ be a nonempty closed convex subset of $E$, and let $T$ be a mapping from $C$ into itself. We denote by $F(T)$ the set of fixed points of $T$. A point $p \in C$ is said to be an asymptotic fixed point of $T$ \cite{Reich} if $C$ contains a sequence $\{x_n\}_{n=1}^\infty$ which converges weakly to $p$ such that $\displaystyle \lim_{n \ra \infty} \lt x_n - Tx_n\rt = 0.$ The set of asymptotic fixed points of $T$ will be denoted by $\tilde{F}(T).$ The asymptotic behavior of relatively nonexpansive mapping was studied in \cite{butnariu1, censor}. A mapping $T$ from $C$ into itself is said to be relatively nonexpansive \cite{su, zegeye1} if the following conditions are satisfied:
\begin{itemize}
\item[(R1)] $F(T) \neq \emptyset$;
\item[(R2)] $\phi(p, Tx) \leq \phi(p, x), \; \forall \; x \in C, \; p \in F(T)$; 
\item[(R3)] $F(T) = \tilde{F}(T)$.
\end{itemize}

\noindent If $E$ is a smooth, strictly convex and reflexive real Banach space, and $A \su E\times E^*$ is a continuous monotone mapping with $A^{-1}(0) \neq \emptyset$ then it proved in \cite{Kohsaka1} that $J_r := (J + rA)^{-1}J,$ for $r>0$ is relatively nonexpansive. Moreover, if $C$ is nonempty closed convex subset of a smooth, strictly convex, and reflexive Banach $E$, and $T :C \ra C$ is relatively nonexpansive then $F(T)$ is \textit{closed and convex} (see, \cite{matsushita1}).

\noindent Let $f: C\times C \ra \R$ be a bifunction. The equilibrium problem for $f$ is to find $x^* \in C$ such that
\begin{equation}
f(x^*, y) \geq 0, \; \forall \; y \in C.
\end{equation}
The set of solutions of (7) is denoted by $EP(f).$ For solving the equilibrium problem for a bifunction $f: C\times C \ra \R,$ let us assume that $f$ satisfies the following conditions:
\begin{itemize}
	\item[(A1)] $f(x, x) = 0$ for all $x \in C$,
	\item[(A2)] $f$ is monotone, that is, $f(x, y) + f(y, x) \leq 0$ for all $x, y \in C$,
	\item[(A3)] for each $x, y, z \in C, \; \displaystyle\lim_{t \ra 0}f(tz+(1-t)x, y) \leq f(x, y)$,
	\item[(A4)] for each $x \in C, y \mapsto f(x, y)$ is convex and lower semicontinuous.
\end{itemize}

\noindent Many authors studied the problem of finding a common element of the set of fixed points of nonexpansive and / or relatively nonexpansive mappings and the set of solutions of an equilibrium problem in the frame work of Hilbert spaces and Banach spaces respectively: see, for instance, \cite{qin, Reich, su2008monotone, zembayashi, zegeye1, ofoedu6, osilike} and the references therein. For finding an element of $F(S)\cap VI(A, C)$, Takahashi and Toyoda \cite{toyoda} introduced the following iterative scheme:
\begin{equation}
x_{n+1} = \al_nx_n+(1-\al_n)SP_C(x_n-\la_nAx_n), \; n \geq 1,
\end{equation}
where $x_0 \in C, \; P_C$ is a metric projection of $H$ onto $C$, $\{\al_n\}_{n \geq 1}$ is a sequence in $(0, 1)$ and $\{\la_n\}$ is a sequence in $(0, 2\al)$, where $\al$ is strong monotonicity constant of $A$.

\noindent Recently, Iiduka and Takahashi \cite{Iiduka} introduced the following iterative scheme:
\begin{equation}
x_{n+1} = \al_nu+\beta_nx_n+\ga_nSP_C(x_n-\la_nAx_n), \; n \geq 1,
\end{equation}
where $u, x_0 \in C,$ and proved the strong convergence theorems for iterative scheme (9) under some conditions on parameters. Furthermore, Tada and Takahashi \cite{tada} introduced the Mann type iterative algorithm for finding a common element of the $EP(f)$ and the set of the common fixed points of nonexpansive mapping and obtained the weak convergence of the Mann type iterative algorithm. In 2007, Yao, Liou and Yao \cite{yao} introduced an iterative process for finding a common element of the set of solutions of the $EP(f)$ and the set of common fixed points of infinitely many nonexpansive mappings in Hilbert spaces. They proved a strong convergence theorem under mild conditions on iterative parameters. Very recently, Moudafi \cite{moudafi1} proposed an iterative algorithm for finding a common element of $GEP(f, B)\cap F(S),$ where $B : C \ra H$ is an $\al$-inverse strongly monotone mapping, and obtained a weak convergence theorem.

\noindent
In this paper, it is our aim to introduce an iterative process for finding a common element of sets of solutions of finite family of generalized equilibrium problem, sets of fixed points of finite family of continuous relatively nonexpansive mappings and sets of zeros of finite family of $\ga$-inverse strongly monotone mappings in Banach spaces.

\subsection*{\normalsize PRELIMINARIES}
\noindent
Let $E$ be a normed linear space with dim $E\geq 2.$ The \emph{modulus of smoothness} of $E$ is the function $\rho_E : [0, \infty) \ra [0, \infty)$ defined by $$\rho_E(\tau) := \sup\left\{\frac{\lt x+y\rt + \lt x-y\rt}{2} - 1: \lt x\rt = 1; \lt y\rt = \tau\right\}.$$
The space $E$ is said to be smooth if $\rho_E(\tau) > 0, \; \forall \; \tau > 0$ and $E$ is called \emph{uniformly smooth} if and only if $\displaystyle\lim_{\tau \ra 0^+}\frac{\rho_E(\tau)}{\tau} = 0.$

\noindent The \emph{modulus of convexity} of $E$ is the function $\delta_E : (0,2] \ra [0,1]$ defined by $$\delta_E(\epsilon) := \inf\left\{1-\lt\frac{x+y}{2}\rt: \lt x\rt = \lt y\rt = 1; \epsilon = \lt x-y\rt\right\}.$$
$E$ is called \emph{uniformly convex} if and only if $\delta_E(\epsilon) > 0$ for every $\epsilon \in (0,2].$ Let $p>1.$ Then $E$ is said to be $p$-\emph{uniformly convex} if there exists a constant $c>0$ such that $\delta(\epsilon) \geq c\epsilon^p$ for all $\epsilon \in [0,2].$ Observe that every p-uniformly convex Banach space is uniformly convex real Banach space.

\noindent It is well known (see for example \cite{Xu2}) that
\[
L_p(l_p) \; \mbox{or} \; W^p_m \; \mbox{is} \; \left\{
\begin{array}{ll}
p-uniformly \; convex & \quad if \; p \geq 2;\\
2-uniformly \; convex & \quad if \; 1<p\leq 2.
\end{array} \right.
\]
In the sequel, we shall make use of the following definitions and lemmas.

\begin{lem} \cite{Xu2}
If $E$ is a uniformly convex real Banach space, then there exists a continuous, strictly increasing and convex function $g : [0, \infty) \ra [0, \infty)$, $g(0) = 0$, such that for all $x, y \in B_r(0) := \{x \in E: \lt x\rt \leq r\}$ and for any $\al \in [0,1]$, we have $$\lt \al x+(1-\al)y\rt^2 \leq \al\lt x\rt^2+(1-\al)\lt y\rt^2-\al(1-\al)g(\lt x-y\rt).$$
\end{lem}

\begin{lem} \cite{Xu2}
Let $E$ be a $2$-uniformly convex real Banach space. Then for all $x, y \in E,$ we have
\begin{equation}
\lt x-y\rt \leq \frac{2}{c^2}\lt Jx-Jy\rt.
\end{equation}
where $J$ is the normalized duality mapping of $E$ and $0<c<1.$
\end{lem}

\begin{lem} \cite{alber}
Let $C$ be a nonempty closed and convex subset of a reflexive, strictly convex, and smooth real Banach space $E$ and let $x \in E.$ Then $\forall \; y \in C,$ $$\phi(y, \Pi_Cx)+\phi(\Pi_Cx, x) \leq \phi(y, x).$$
\end{lem}

\begin{lem} \cite{kamimura}
Let $C$ be a nonempty closed convex subset of a smooth, uniformly convex Banach space $E$. Let $\{x_n\}_{n=1}^\infty$ and $\{y_n\}_{n=1}^\infty$ be sequences in $E$ such that either $\{x_n\}_{n=1}^\infty$ or $\{y_n\}_{n=1}^\infty$ is bounded. If $\displaystyle\lim_{n \ra \infty}\phi(x_n, y_n) = 0,$ then $\displaystyle\lim_{n \ra \infty}\lt x_n - y_n\rt = 0.$
\end{lem}

\begin{lem} \cite{alber}
Let $C$ be a convex subset of a smooth real Banach space $E$. Let $x \in E$. Then $x_0 = \Pi_Cx$ if and only if $$\ltl z-x_0, Jx_0-Jx\rtl \geq 0, \forall z \in C.$$
\end{lem}

\noindent
We denote by $N_C(v)$ the \textit{normal cone} for $C$ at a point $v \in C$, that is $N_C(v) := \{x^* \in E^* : \ltl v-y, x^* \rtl \geq 0 \; \mbox{for all} \; y \in C\}$. In the sequel we shall use the following lemma.

\begin{lem} (Rockafellar \cite{Rockafellar1})
Let $C$ be a nonempty closed convex subset of a real Banach space $E$ and let $A$ be a monotone and hemicontinuous opeartor of $C$ into $E^*$ with $C = D(A).$ let $B \su E\times E^*$ be an operator defined as follows:
\begin{equation}
Bv := \left\{
\begin{array}{ll}
Av+N_C(v), & v \in C,\\
\emptyset, & v \in C.
\end{array}
\right.
\end{equation}
Then $B$ is maximal monotone and $B^{-1}(0) = VI(A, C)$.
\end{lem}

\noindent We make use of the function $V: E\times E^* \ra \R$ defined by $$V(x, x^*) = \lt x\rt^2-2\ltl x, x^*\rtl+\lt x^*\rt^2, \; \mbox{for all} \; x \in E \;\; \mbox{and} \;\; x^* \in E^*,$$
studied by Alber \cite{alber}. That is $V(x, x^*) = \phi(x, J^{-1}x^*)$ for all $x \in E$ and $x^* \in E^*.$ We know the following Lemma.

\begin{lem} \cite{alber}
Let $E$ be a reflexive strictly convex and smooth real Banach space with $E^*$ as its dual. Then $$V(x, x^*) + 2\ltl j^{-1}x^*-x, y^*\rtl \leq V(x, x^*+y^*),$$ for all $x \in E$ and $x^*, y^* \in E^*.$
\end{lem}

\begin{lem} \cite{blum}
Let $C$ be closed convex subset of a smooth srictly convex and reflexive real Banach space $E$. Let $f$ be a bifunction from $C\times C$to $\R$ satisfying (A1)-(A4). Then, for $r > 0$ and $x \in E,$ there exists $z \in C$ such that 
\begin{equation}
f(z, y)+\frac{1}{r}\ltl y-z, Jz-Jx\rtl \geq 0, \; \forall \; y \in C.
\end{equation}
\end{lem}

\begin{re}
Replacing $x$ with $J^{-1}(Jx-rB(x))$ in (12), where $B$ is monotone mapping from $C$ into $E^*$, then there exists $z \in C$ such that $$f(z, y)+\ltl Bx, y-z\rtl+\frac{1}{r}\ltl y-z, Jz-Jx\rtl \geq 0, \; \forall \; y \in C.
$$
\end{re}

\noindent By a similar argument of the proof of Lemma 2.8 and Remark 2.9 of \cite{takahashi2008strong}, we have the following Lemmas.

\begin{lem} \cite{zembayashi}
Let $C$ be a nonempty closed convex subset of a smooth, strictly convex and reflexive real Banach space $E$. Let $F$ be a bifunction from $C\times C$ to $\R$ satisfying (A1)-(A4) and $B : C \ra E^*$ be a monotone mapping. For $r > 0$ and $x \in E,$ define a mapping $T_r : E \ra C$ as follows:
$$T_rx := \{z \in C: f(z, y)+\ltl Bx, y-z\rtl+\frac{1}{r}\ltl y-z, Jz-Jx\rtl \geq 0, \; \forall \; y \in C\}$$
for all $x \in E.$ Then the following hold:
\begin{itemize}
	\item[(1)] $T_r$ is single valued
	\item[(2)] $T_r$ is firmly nonexpansive type mapping, i.e., for all $x, y \in E$,
	$$\ltl T_rx-T_ry, JT_rx - JT_ry\rtl \leq \ltl T_rx - T_ry, Jx-Jy\rtl;$$
	\item[(3)] $F(T_r) = GEP(F, B)$;
	\item[(4)] $GEP(F, B)$ is closed and convex.
\end{itemize}
\end{lem}

\begin{lem} \cite{zembayashi}
Let $C$ be a nonempty closed convex subset of a smooth, strictly convex and reflexive Banach space $E$. Let $F$ be a bifunction from $C\times C$ to $\R$ satisfying (A1)-(A4). For $r>0,$ $x \in E$ and $p \in F(T_r)$, we have that $$\phi(q, T_rx)+\phi(T_rx, x) \leq \phi(q, x).$$
\end{lem}

\subsection*{\normalsize RESULTS}
\noindent
Let $C$ be a nonempty closed convex subset of $2$-uniformly convex and uniformly smooth real Banach space $E$. Let $f_1, f_2 : C\times C \ra \R, \; k = 1,2, \cdots, q$ be bifunctions and $B_1, B_2 : C \ra E^*$ be monotone mappings. Let $T_j : C \ra C, \; j = 1, 2, \cdots,d$ be finite family of continuous relatively nonexpansive mappings and $A_i : C \ra E^*, \; i = 1,2, \cdots, m$ be finite family of $\ga_i$-inverse strongly monotone operators with constants $\ga_i \in (0,1), \; i = 1, 2, \cdots, m;$ then in what follows, we shall study the following iteration process.

\begin{equation}
\left\{ \begin{array}{ll}
x_0 \in C_0 = C, \quad \mbox{chosen arbitrarily},\\
z_n = \Pi_{C}J^{-1}(Jx_n - \la_nA_{n+1}x_n),\\
y_n = J^{-1}(\alpha_{n}Jx_n+(1-\al_n)JT_{n+1}z_n); u_n, v_n \in C \; \mbox{s.t.}\\
f_1(u_n,y)+\ltl B_1y_n, y-u_n\rtl+\frac{1}{r_n}\ltl y-u_n, Ju_n-Jy_n\rtl \geq 0
, \forall \; y \in C,\\
f_2(v_n,y)+\ltl B_2y_n, y-v_n\rtl+\frac{1}{r_n}\ltl y-v_n, Jv_n-Jy_n\rtl \geq 0
, \forall \; y \in C,\\
w_n = J^{-1}(\beta Ju_n+(1-\beta) Jv_n),\\
C_{n+1} = \{z \in C_n: \; \phi(z, w_n) \leq \phi(z, x_n)\},\\
x_{n+1} = \Pi_{C_{n+1}}(x_0), \; n\geq 0,
\end{array}
\right.
\end{equation}

\noindent
where $A_n = A_{n(\bmod m)}$, $T_n = T_{n(\bmod m)}$ and $J$ is the normalized duality mapping on $E$; $\{r_n\}_{n \geq 1} \su [c_1, \infty)$ for some $c_1 > 0$, $\beta$, $\al_n \in (0,1)$ for all $n \in \N$ such that $\ds\liminf_{n \ra \infty}\al_n(1-\al_n) > 0$; and $\{\la_n\}_{n\geq 1}$ is a sequence in $[a, b]$ for some $0<a<b<\frac{c^2\ga}{2}$, where $c$ is the $2$-uniformly convex constant of $E$ and $\ga = \ds\min_{1\leq i\leq m}\ga_i$.

\noindent
We shall define $$T_{k,r}x := \{z \in C: f_k(z,y)+\ltl B_kx, y-z\rtl+\frac{1}{r}\ltl y-z, Jz-Jx\rtl \geq 0, \; \forall \; y \in C\}$$ for all $x\in E, k = 1, 2.$

\begin{lem}
Let $C$ be a nonempty closed convex subset of $2$-uniformly convex and uniformly smooth real Banach space $E$. Let $f_1, f_2 : C \times C \ra \R$ be bifunctions satisfying (A1) - (A4) and $B_1, B_2: C \ra E^*$ be continuous monotone mappings. Let $T_j : C \ra C, \; j = 1, 2, \ldots, d$ be a finite family of relatively nonexpansive mappings and $A_i : C \ra E^*, i = 1, 2, \ldots, m$ be a finite family of $\ga_i$-inverse strongly monotone operators with constants $\ga_i \in (0,1), \; i = 1, 2, \ldots, m$. Let $F := \ds\left[\bigcap_{j=1}^dF(T_j)\right]\cap\left[\bigcap_{i=1}^mA_i^{-1}(0)\right]\cap\left[\bigcap_{k=1}^2GEP(f_k, B_k)\right] \neq \emptyset$ and let $\{x_n\}$ be a sequence defined by (13). Then the sequence $\{x_n\}$ is well defined for each $n \geq 0$.
\end{lem}

\begin{prf}
We first show that $C_n$ is closed and convex for all $n\geq 1.$ It is obvious that $C_0 = C$ is closed and convex from the definition. Suppose that $C_n$ is closed convex for some $n \geq 1.$ From the definition of $C_{n+1},$ we have that $z \in C_{n+1}$ implies $\phi(z, w_n)\leq \phi(z, x_n).$ This is equivalent to $$\lt z\rt^2-2\ltl z, Jw_n\rtl+\lt w_n\rt^2\leq\lt z\rt^2-2\ltl z, Jx_n\rtl+\lt x_n\rt^2$$ $$ \mbox{Which gives} \quad 2(\ltl z, Jx_n\rtl -\ltl z, Jw_n\rtl) = 2(\ltl z, Jx_n - Jw_n\rtl) - \lt x_n\rt^2-\lt w_n\rt^2 \leq 0$$ It follows that $C_{n+1}$ is closed and convex, hence $\Pi_{C_{n+1}}$ is well defined for all $n\geq 0.$

\noindent
Next we prove that $F \su C_n$ for all $n \geq 0.$ From the assumption, we see that $F \su C_0 = C.$ Suppose that $F \su C_n$ for some $n \geq 1.$ Now, for $p \in F,$ the property of $G$ and Lemma 2 give that
\begin{equation}
\begin{aligned}
\phi(p, w_n) =& \phi(p, J^{-1}(\beta Ju_n + (1-\beta)Jv_n))\\
=& \lt p\rt^2-2\ltl p, \beta Ju_n + (1-\beta)Jv_n\rtl+\lt\beta Ju_n + (1-\beta)Jv_n\rt^2\\
\leq& \lt p\rt^2-2\beta\ltl p, Ju_n\rtl-2(1-\beta)\ltl p, Jv_n\rtl+\beta\lt Ju_n\rt^2+(1-\beta)\lt Jv_n\rt^2\\
=& \beta\phi(p, u_n) + (1-\beta)\phi(p, v_n)\\ 
=& \beta\phi(p, T_{1,r_n}y_n) + (1-\beta)\phi(p, T_{2,r_n}y_n)\\
\leq&  \beta\phi(p, y_n) + (1-\beta)\phi(p, y_n) = \phi(p, y_n)\\
=& \phi(p, J^{-1}(\al_nJx_n+(1-\al_n)JT_{n+1}z_n))\\
=& \lt p\rt^2-2\ltl p, \al_nJx_n+(1-\al_n)JT_{n+1}z_n\rtl+\lt \al_nJx_n+(1-\al_n)JT_{n+1}z_n\rt^2\\
\leq& \lt p\rt^2-2\al_n\ltl p, Jx_n\rtl-2(1-\al_n)\ltl p, JT_{n+1}z_n\rtl+\al_n\lt Jx_n\rt^2\\
+& (1-\al_n)\lt JT_{n+1}z_n\rt^2\\
=& \al_n\phi(p, x_n)+(1-\al_n)\phi(p, T_{n+1}z_n)\\
\leq& \al_n\phi(p, x_n)+(1-\al_n)\phi(p, z_n).
\end{aligned}
\end{equation}
Moreover, by Lemma 4 and Lemma 8 we get that 
\begin{equation}
\begin{aligned}
\phi(p, z_n) =& \phi(p, \Pi_CJ^{-1}(Jx_n - \la_nA_{n+1}x_n))\\
\leq& \phi(p, J^{-1}(Jx_n - \la_nA_{n+1}x_n))\\
=& V(p, Jx_n - \la_nA_{n+1}x_n)\\
\leq& V(p, (Jx_n-\la_nA_{n+1})+\la_nA_{n+1}x_n)-2\ltl J^{-1}(Jx_n-\la_nA_{n+1}x_n)-p, \la_nA_{n+1}x_n\rtl\\
=& V(p, Jx_n) - 2\la_n\ltl J^{-1}(Jx_n - \la_nA_{n+1}x_n)-p, A_{n+1}x_n\rtl\\
=& \phi(p, Jx_n) - 2\la_n\ltl x_n-p, A_{n+1}x_n\rtl-2\la_n\ltl J^{-1}(Jx_n - \la_nA_{n+1}x_n)-x_n, A_{n+1}x_n\rtl\\
\leq& \phi(p, x_n) - 2\la_n\ltl x_n-p, A_{n+1}\rtl+2\ltl J^{-1}(Jx_n-\la_nA_{n+1}x_n)-x_n, -\la_nA_{n+1}x_n\rtl
\end{aligned}
\end{equation}
Thus, since $p \in \ds\cap_{i=1}^mA_i^{-1}(0)$ and $\ga = \min\ga_i$, we have that 
$$\ltl x-p, A_ix\rtl \geq \ga\lt A_i x\rt^2 \quad \mbox{for} \quad i = 1, 2, \cdots, m.$$
Thus, we have from (15) that
\begin{equation}
\phi(p, z_n) \leq \phi(p, x_n)-2\la_n\ga\lt A_{n+1}x_n\rt^2+2\ltl J^{-1}(Jx_n - \la_nA_{n+1}x_n)-x_n, -\la_nA_{n+1}x_n\rtl
\end{equation}
Therefore, from (10), (16) and $\la_n \leq \frac{c^2\ga}{2}$ we obtain that
\begin{equation}
\begin{aligned}
\phi(p, z_n) \leq& \phi(p, x_n)-2\la_n\ga\lt A_{n+1}x_n\rt^2+\frac{4\la_n^2}{c^2}\lt A_{n+1}x_n\rt^2\\
=& \phi(p, x_n)+2\la_n\left(\frac{2}{c^2}\la_n-\ga\right)\lt A_{n+1}x_n\rt^2 \leq \phi(p, x_n).
\end{aligned}
\end{equation}
Substituting (17) into (14), we have $$\phi(p, w_n) \leq \phi(p, x_n),$$
that is $p \in C_{n+1}$. This implies, by induction, that $F \su C_n$ and the sequence $\{x_n\}_{n=0}^\infty$ generated by (13) is well defined for all $n \geq 0.$
\end{prf}

\begin{theo}
Let $C$ be a nonempty closed convex subset of $2$-uniformly convex and uniformly smooth real Banach space $E$. Let $f_1, f_2 : C \times C \ra \R$ be bifunctions satisfying (A1) - (A4) and $B_1, B_2 : C \ra E^*$ be continuous monotone mappings. Let $T_j : C \ra C, \; j = 1, 2, \ldots, d$ be a finite family of continuous relatively nonexpansive mappings and $A_i : C \ra E^*, \; i = 1, 2, \ldots, m$ be a finite family of $\ga_i$-inverse strongly monotone operators with constants $\ga_i \in (0,1), \; i = 1, 2, \ldots, m.$ Let $F := \ds\left[\bigcap_{j=1}^dF(T_j)\right]\cap\left[\bigcap_{i=1}^mA_i^{-1}(0)\right]\cap\left[\bigcap_{k=1}^2GEP(f_k, B_k)\right] \neq \emptyset$. Let $\{x_n\}_{n\geq 0}$ be a sequence defined by (13). Then, the sequence $\{x_n\}_{n\geq 0}$ converges to some element of $F$.
\end{theo}

\begin{prf}
We have from Lemma 12 that $F \su C_n, \; \forall \; n \geq 0$ and $x_n$ is well defined for each $n \geq 0$. From $x_n = \Pi_{C_n}(x_0)$, and Lemma 4, we have $$\phi(p, x_0) = \phi(\Pi_{C_n}x_0, x_0) \leq \phi(p, x_0) - \phi(p, x_n) \leq \phi(p, x_0),$$ for each $p \in F \su C_n$ and $n \geq 0$. Thus the sequence $\{\phi(x_n, x_0)\}_{n=0}^\infty$ is bounded. 

\noindent
Furthermore, since $x_n = \Pi_{C_n}(x_0)$ and $x_{n+1} = \Pi_{C_{n+1}}(x_0) \in C_{n+1} \su C_n$ we have that $$\phi(x_n, x_0) \leq \phi(x_{n+1}, x_0), \; \forall \; n\geq 0,$$ which implies that $\{\phi(x_n, x_0)\}_{n\geq 0}$ is increasing and hence $\ds \lim_{n\ra \infty}\phi(x_n, x_0)$ exists. Similarly, by Lemma 4, we have, for any positive integer $l$, that
\[
\begin{aligned}
\phi(x_{n+l}, x_n) =& \phi(x_{n+l}, \Pi_{C_n}x_0)\\
\leq& \phi(x_{n+l}, x_0) - \phi(\Pi_{C_n}x_0, x_0)\\
=& \phi(x_{n+l}, x_0) - \phi(x_n, x_0) \; \forall \; n\geq 0.
\end{aligned}
\]

\noindent
Since $\ds\lim_{n\ra\infty}\phi(x_n, x_0)$ exists, we have that $\ds\lim_{n\ra\infty}\phi(x_{n+l},x_n) = 0$. Thus, Lemma 5 implies that
\begin{equation}
\lim_{n\ra\infty}\lt x_{n+l} - x_n\rt = 0,
\end{equation}
and hence $\{x_n\}$ is cauchy. Therefore, there exists a point $x^* \in C$ such that $x_n \ra x^*$ as $n \ra \infty$. Since $x_{n+1} = \Pi_{C_{n+1}}x_0 \in C_{n+1}$, we have $\phi(x_{n+1}, w_n)\leq \phi(x_{n+1}, x_n), \; \forall \; n\geq 0.$\\ Thus, by (18) and Lemma 5 we get that 
\begin{equation}
\lim_{n\ra\infty}\lt x_{n+1} - w_n\rt = 0,
\end{equation}
and hence $\lt x_n - w_n\rt \leq \lt x_n - x_{n+1}\rt + \lt x_{n+1} - w_n\rt \ra 0$ as $n \ra \infty$, which implies that $w_n \ra x^*$ as $n \ra \infty$. Furthermore, the uniform continuity of $J$ on bounded sets, gives that
\begin{equation}
\lim_{n\ra\infty}\lt Jx_{n+l} - Jw_n\rt = 0,
\end{equation}

\noindent
We note that if $E$ is uniformly smooth then $E^*$ is uniformly convex. Thus, using property of $\phi$ and Lemma 2 we have, for all $p \in F$, that
\begin{equation}
\begin{aligned}
\phi(p, y_n) =& \phi(p, J^{-1}(\al_nJx_n+(1-\al_n)JT_{n+1}z_n)\\
=& \lt p\rt^2-2\ltl p, \al_nJx_n+(1-\al_n)JT_{n+1}z_n\rtl+\lt \al_nJx_n+(1-\al_n)JT_{n+1}z_n\rt^2\\
\leq& \lt p\rt^2-2\al_{n}\ltl p, Jx_n\rtl-2(1-\al_{n})\ltl p, JT_{n+1}z_n\rtl+\al_{n}\lt Jx_n\rt^2+(1-\al_n)\lt JT_{n+1}z_n\rt^2\\
-& \al_n(1-\al_n)g(\lt Jx_n-JT_{n+1}z_n\rt)\\
=& \al_n\phi(p, x_n)+(1-\al_n)\phi(p, T_{n+1}z_n)-\al_n(1-\al_n)g(\lt Jx_n-JT_{n+1}z_n\rt)\\
\leq& \al_n\phi(p, x_n)+(1-\al_n)\phi(p, z_n)-\al_n(1-\al_n)g(\lt Jx_n-JT_{n+1}z_n\rt),
\end{aligned}
\end{equation}
Thus, from (17) and (21) we have that 
\begin{equation}
\phi(p, y_n) \leq \phi(p, x_n)+2(1-\al_n)\la_n\left(\frac{2}{c^2}\la_n-\ga\right)\lt A_{n+1}\rt^2 - \al_n(1-\al_n)g(\lt Jx_n-JT_{n+1}z_n\rt)
\end{equation}
On the other hand from Lemma 11 we get that
\begin{equation}
\begin{aligned}
\phi(p, w_n) =& \phi(p, J^{-1}(\beta Ju_n + (1-\beta)Ju_n))\\
\leq& \beta\phi(p, u_n) + (1-\beta)\phi(p, v_n)\\
=& \beta\phi(p, T_{1,r}y_n) + (1-\beta)\phi(p, T_{2,r}y_n) \leq \phi(p, y_n).
\end{aligned}
\end{equation}
Substituting (22) into (23) we get that
\begin{equation}
\phi(p, w_n) \leq \phi(p, x_n)+2(1-\al_n)\la_n\left(\frac{2}{c^2}\la_n-\ga\right)\lt A_{n+1}\rt^2-\al_n(1-\al_n)g(\lt Jx_n-JT_{n+1}z_n\rt).
\end{equation}
Now, using the fact that $\la_n < \frac{c^2\ga}{2},$ the inequality (24) implies that
\begin{equation}
\begin{aligned}
\al_n(1-\al_n)g(\lt Jx_n-JT_{n+1}z_n\rt) \leq& \phi(p, x_n)-\phi(p, w_n)\\
=& \lt p\rt^2-2\ltl p, Jx_n\rtl+\lt x_n\rt^2-\lt p\rt^2+2\ltl p, Jw_n\rtl-\lt w_n\rt^2\\
=& \lt x_n\rt^2 - \lt w_n\rt^2-2\ltl p, Jx_n - Jw_n\rtl\\
\leq& \lt x_n - w_n\rt(\lt x_n\rt+\lt w_n\rt)+2\lt p\rt\lt Jx_n-Jw_n\rt\\
\leq& M_0(\lt x_n-w_n\rt+\lt Jx_n-Jw_n\rt),
\end{aligned}
\end{equation}
for some $M_0 > 0.$ Thus, since $\ds\lim_{n \ra \infty}\lt x_n-w_n\rt = 0, \; \ds\lim_{n \ra \infty}\lt Jx_n-Jw_n\rt=0$ and we obtain 
\begin{equation}
\phi(p, x_n)-\phi(p, w_n) \lra 0, \; \mbox{as} \; n \lra \infty,
\end{equation}
and hence inequality (25) implies that $g(\lt Jx_n-JT_{n+1}z_n\rt) \ra 0 \; \mbox{as} \; n \ra \infty.$ Therefore, from the property of $g$ we get that $\lt Jx_n - JT_{n+1}z_n\rt \ra 0 \; \mbox{as} \; n \ra \infty.$
Furthermore, since $J^{-1}$ is also uniformly norm-to-norm continuous on bounded sets, we see that
\begin{equation}
\lim_{n \ra \infty}\lt x_n-T_{n+1}z_n\rt = 0.
\end{equation}
Moreover, from (24) we have that $$(1-\al_n)2\la_n\left(\ga-\frac{2}{c^2}\la_n\right)\lt A_{n+1}x_n\rt^2 \leq \phi(p, x_n)-\phi(p, w_n)$$
which yields that
\begin{equation}
\lim_{n\ra \infty}\lt A_{n+1}x_n\rt = 0
\end{equation}
Now, Lemma 4, Lemma 8 and (17) imply that
\[
\begin{aligned}
\phi(x_n, z_n) =& \phi(x_n, \Pi_CJ^{-1}(Jx_n - \la_nA_{n+1}x_n))\\
\leq& \phi(x_n, J^{-1}(Jx_n-\la_nA_{n+1}x_n))\\
=& V(x_n, Jx_n - \la_nA_{n+1}x_n)\\
=& V(x_n, (Jx_n-\la_nA_{n+1}x_n)+\la_nA_{n+1}x_n)-2\ltl J^{-1}(Jx_n-\la_nA_{n+1}x_n)-x_n, \la_nA_{n+1}x_n\rtl\\
=& \phi(x_n, x_n)+2\ltl J^{-1}(Jx_n-\la_nA_{n+1}x_n)-x_n, -\la_nA_{n+1}x_n\rtl\\
=& 2\ltl J^{-1}(Jx_n-\la_nA_{n+1}x_n)-x_n, -\la_nA_{n+1}x_n\rtl\\
\leq& 2\lt J^{-1}(Jx_n-\la_nA_{n+1})-J^{-1}Jx_n\rt\lt\la_nA_{n+1}\rt\\
=& 2\lt J^{-1}\la_nA_{n+1}x_n\rt\lt\la_nA_{n+1}x_n\rt\\
\leq& \frac{4b^2}{c^2}\lt A_{n+1}x_n\rt^2.
\end{aligned}
\]
that is
\begin{equation}
\phi(x_n, z_n) \leq \frac{4b^2}{c^2}\lt A_{n+1}x_n\rt^2.
\end{equation}
It follows from (28), (29) and Lemma 5 that
\begin{equation}
\lim_{n \ra \infty}\lt x_n - z_n\rt = 0;
\end{equation}
and hence $z_n \ra p$ as $n \ra \infty.$

\noindent
We now show that $p \in \ds\cap_{j=1}^dF(T_i).$ Observe that from (27) and (30) we obtain that
$$\lt T_{n+1}z_n-z_n\rt \leq \lt T_{n+1}z_n - x_n\rt+\lt z_n-x_n\rt \ra 0 \; \mbox{as} \; n \ra \infty.$$
Hence,
\begin{equation}
\lim_{n \ra \infty}T_{n+1}z_n = p.
\end{equation}
Let $\{n_l\}_{l \geq 1}\su\N$ be such tat $T_{n_l+1} = T_1$ for all $l \in \N$, then since $z_{n_l} \ra p$ as $l \ra \infty$, we obtain from (31), using continuity of $T_1$, that
$$p=\lim_{l \ra \infty}T_{n_l+1}z_{n_l} = \lim_{l \ra \infty}T_1z_{n_l} = T_1p.$$
Similarly, if $\{n_k\}_{k\geq 1} \su \N$ is such that $T_{n_k+1} = T_2$ for all $k \in \N$, then we have again that $$p = \lim_{k\ra\infty}T_{n_k+1}z_{n_k} = \lim_{k\ra\infty}T_2z_{n_k} = T_2x^*.$$ Continuing, we obtain that $T_jp = p, \; j = 1, 2, \cdots, d.$ Hence, $p \in \ds\bigcap_{j=1}^dF(T_j).$

\noindent
Next we show that $p \in \ds\bigcap_{i=1}^mA_i^{-1}(0).$ Since $A_i$ is $\ga$-inverse strongly monotone for $i = 1, 2, \cdots, m,$ we have that $A_i$, is $\frac{1}{\ga}$-Lipschitz continuous. Thus,
\begin{equation}
\lt A_{n+1}x_n-A_{n+1}p\rt \leq \frac{1}{\ga}\lt x_n - p\rt \ra 0, \; \mbox{as} \; n \ra \infty.
\end{equation}
Hence, from (32) and (28), we obtain that $$\lt A_{n+1}p\rt \leq \lt A_{n+1}x_n-A_{n+1}p\rt + \lt A_{n+1}x_n\rt \ra 0 \; \mbox{as} \; n \ra \infty.$$
Consequently, we get that
$$\lim_{n \ra \infty}A_{n+1}p = 0.$$
Let $\{n_s\}_{s\geq 1} \su \N$ be such that $A_{n_s+1} = A_1$ for all $s \in \N.$ Then, $A_1p = \ds\lim_{s \ra \infty}A_{n_s+1}p = 0.$\\
Similarly, we have that $A_ip = 0$ for $i = 2, \cdots, m.$ Thus, $p \in \ds\bigcap_{i=1}^mA_{i}^{-1}(0).$

\noindent
Furthermore, we show that $p \in \ds\bigcap GEP(f_k, B_k) = F(T_{k,r}), \; k = 1, 2.$ Let $p \in F.$ From $u_n = T_{1,r_n}y_n,$ Lemma 11, (21) and the fact that $x_n \ra p, \; z_n \ra p$ as $n \ra \infty,$ we obtain
\begin{equation}
\begin{aligned}
\phi(p, u_n) =& \phi(p, T_{1,r_n}y_n) \leq \phi(p, y_n)\\
\leq& \al_n\phi(p, x_n)+(1-\al_n)\phi(p, z_n)\\
\leq& \phi(p, x_n)+\phi(p, z_n) \ra 0 \; \mbox{as} \; n \ra \infty.
\end{aligned}
\end{equation}
Thus, by Lemma 5 and (33), $u_n \ra p$ and $y_n \ra p$ as $n \ra \infty.$ These imply that, $$\lt u_n - y_n\rt \ra 0 \; \mbox{as} \; n \ra \infty.$$
Consequently, $\ds\lim_{n\ra\infty}\lt Ju_n - Jy_n\rt = 0.$ Hence,
\begin{equation}
\lim_{n \ra \infty}\frac{\lt Ju_n - Jy_n\rt}{r_n} = 0.
\end{equation}
But from (A2) we note that $$\ltl B_1y_n, v-u_n\rtl+\frac{1}{r_n}\ltl v-u_n, Ju_n-Jy_n\rtl \geq -f_1(u_n, v) \geq f_1(v, u_n) \; \forall \; v \in C,$$ and hence
\begin{equation}
\ltl B_1y_n, v-u_n\rtl+\ltl v-u_n, \frac{Ju_n-Jy_n}{r_n}\rtl \geq f_1(v, u_n) \; \forall \; v \in C.
\end{equation}
Put $z_t = tv + (1-t)p$ for all $t \in (0, 1]$ and $v \in C.$\\
Consequently, we get that $z_t \in C.$ From (35), it follows that
\[
\begin{aligned}
\ltl B_1z_t, z_t-u_n\rtl -& \ltl B_1z_t, z_t-u_n\rtl \geq -\ltl B_1y_n, z_t-u_n\rtl-\ltl z_t-u_n, \frac{Ju_n-Jy_n}{r_n}\rtl+f_1(z_t, u_n)\\
\mbox{This implies that}\\
\ltl B_1z_t, z_t-u_n\rtl \geq& \ltl B_1z_t, z_t-u_n\rtl -\ltl B_1y_n, z_t-u_n\rtl-\ltl z_t-u_n, \frac{Ju_n-Jy_n}{r_n}\rtl+f_1(z_t, u_n)\\
=& \ltl B_1z_t-B_1u_n, z_t-u_n\rtl + \ltl B_1u_n-B_1y_n, z_t-u_n\rtl\\
-&\ltl z_t-u_n, \frac{Ju_n-Jy_n}{r_n}\rtl+f_1(z_t, u_n).
\end{aligned}
\]
By the continuity of $B_1$ and the fact that $u_n \ra p, y_n \ra p$ as $n \ra \infty,$ we obtain that
\begin{equation}
B_1u_n - B_1y_n \ra 0 \;\; \mbox{as} \;\; n \ra \infty.
\end{equation}
Since $B_1$ is monotone, we have that $\ltl B_1z_t - B_1u_n, z_t - u_n\rtl \geq 0.$ Using this, (34) and (36), it follows from (A4) and (36) that $$f_1(z_t, p) \leq \liminf_{n \ra \infty}f_1(z_t, u_n) \leq \lim_{n \ra \infty}\ltl B_1z_t, z_t-u_n\rtl = \ltl B_1z_t, z_t-p\rtl.$$
Now, from (A1) and (A4) we get that
\[
\begin{aligned}
0 =& f_1(z_t, z_t) \leq tf_1(z_t, v)+(1-t)f_1(z_t, p)\\
\leq& tf_1(z_t, v)+(1-t)\ltl B_1z_t, z_t-p\rtl\\
\leq& tf_1(z_t, v)+(1-t)\ltl B_1z_t, tv+(1-t)p-p\rtl\\
=& tf_1(z_t, v)+(1-t)t\ltl B_1z_t, v-p\rtl.
\end{aligned}
\]
and hence $$f_1(z_t, v)+(1-t)\ltl B_1z_t, v-p\rtl \geq 0.$$
Letting $t \ra 0,$ we have
$$f_1(p, v)+\ltl B_1p, v-p\rtl \geq 0.$$
This implies that $p \in GEP(f_1, B_1).$ Similarly, considering $v_n = T_{2,r_n}y_n$, the same argument gives that $p \in GEP(f_2, B_2).$

\noindent
Finally, we prove that $p = \Pi_{F}(x_0).$ From $x_n = \Pi_{C_n}(x_0)$, we have $$\ltl Jx_0-Jx_n, x_n-z\rtl \geq 0 \; \forall \; z \in C_n.$$
Since $F \in C_n,$ we also have that
\begin{equation}
\ltl Jx_0-Jx_n, x_n-p\rtl \geq 0 \; \forall \; p \in F.
\end{equation}
By taking limits in (37), one has $$\ltl Jx_0-Jx^*, x^*-p\rtl \geq 0, \; \forall \; p \in F.$$
Now, by Lemma 6 we have that $x^* = \Pi_{F}x_0$. This completes the proof.
\end{prf}

\noindent
Strong convergence theorem for approximating a common element of sets of solutions of two generalized equilibrium problems and the sets of fixed points of finite family of relatively nonexpansive mappings in Banach spaces may not require that $E$ is a $2$-uniformly convex real Banach space. In fact, we have the following Theorem.

\begin{theo}
Let $C$ be a nonempty closed convex subset of a uniformly convex and uniformly smooth real Banach space $E$. Let $f_1, f_2 : C \times C \ra \R$ be bifunctions satisfying (A1) - (A4) and $B_1, B_2 : C \ra E^*$ be continuous monotone mappings. Let $T_j : C \ra C, \; j = 1, 2, \ldots, d$ be a finite family of continuous relatively nonexpansive mappings. Let $F := \ds\left[\bigcap_{j=1}^dF(T_j)\right]\cap\left[\bigcap_{k=1}^2GEP(f_k, B_k)\right] \neq \emptyset$. Let $\{x_n\}_{n\geq 0}$ be a sequence defined by

\begin{equation}
\left\{ \begin{array}{ll}
x_0 \in C_0 = C, \quad \mbox{chosen arbitrarily},\\
y_n = J^{-1}(\alpha_{n}Jx_n+(1-\al_n)JT_{n+1}x_n); u_n, v_n \in C \; \mbox{s.t.}\\
f_1(u_n,y)+\ltl B_1y_n, y-u_n\rtl+\frac{1}{r_n}\ltl y-u_n, Ju_n-Jy_n\rtl \geq 0
, \forall \; y \in C,\\
f_2(v_n,y)+\ltl B_2y_n, y-v_n\rtl+\frac{1}{r_n}\ltl y-v_n, Jv_n-Jy_n\rtl \geq 0
, \forall \; y \in C,\\
w_n = J^{-1}(\beta Ju_n+(1-\beta) Jv_n),\\
C_{n+1} = \{z \in C_n: \; \phi(z, w_n) \leq \phi(z, x_n)\},\\
x_{n+1} = \Pi_{C_{n+1}}(x_0), \; n\geq 0,
\end{array}
\right.
\end{equation}

\noindent
where $T_n = T_{n(\bmod m)}$ and $J$ is the normalized duality mapping on $E$; $\{r_n\}_{n \geq 1} \su [c_1, \infty)$ for some $c_1 > 0$, $\beta$, $\al_n \in (0,1)$ for all $n \in \N$ such that $\ds\liminf_{n \ra \infty}\al_n(1-\al_n) > 0$. Then, the sequence $\{x_n\}_{n\geq 0}$ converges to some element of $F$.
\end{theo}

\begin{prf}
Put $A_i \equiv 0, \; i = 1,2,\ldots, m$ in Theorem 1. Then, we get that $z_n = x_n;$ and the method of proof of Theorem 1 gives the required assertion without the requirement that $E$ is a $2-$uniformly convex real Banach space.
\end{prf}

\noindent
If, in Theorem 2, we have that $B_1 \equiv B_2 \equiv 0,$then we get the following corollary

\begin{cor}
Let $C$ be a nonempty closed convex subset of a uniformly convex and uniformly smooth real Banach space $E$. Let $f_1, f_2 : C \times C \ra \R$ be bifunctions satisfying (A1) - (A4). Let $T_j : C \ra C, \; j = 1, 2, \ldots, d$ be a finite family of continuous relatively nonexpansive mappings. Let $F := \ds\left[\bigcap_{j=1}^dF(T_j)\right]\cap\left[\bigcap_{k=1}^2EP(f_k)\right] \neq \emptyset$. Let $\{x_n\}_{n\geq 0}$ be a sequence defined by

\begin{equation}
\left\{ \begin{array}{ll}
x_0 \in C_0 = C, \quad \mbox{chosen arbitrarily},\\
y_n = J^{-1}(\alpha_{n}Jx_n+(1-\al_n)JT_{n+1}x_n); u_n, v_n \in C \; \mbox{s.t.}\\
f_1(u_n,y)+\frac{1}{r_n}\ltl y-u_n, Ju_n-Jy_n\rtl \geq 0
, \forall \; y \in C,\\
f_2(v_n,y)+\frac{1}{r_n}\ltl y-v_n, Jv_n-Jy_n\rtl \geq 0
, \forall \; y \in C,\\
w_n = J^{-1}(\beta Ju_n+(1-\beta) Jv_n),\\
C_{n+1} = \{z \in C_n: \; \phi(z, w_n) \leq \phi(z, x_n)\},\\
x_{n+1} = \Pi_{C_{n+1}}(x_0), \; n\geq 0,
\end{array}
\right.
\end{equation}

\noindent
where $T_n = T_{n(\bmod m)}$ and $J$ is the normalized duality mapping on $E$; $\{r_n\}_{n \geq 1} \su [c_1, \infty)$ for some $c_1 > 0$, $\beta$, $\al_n \in (0,1)$ for all $n \in \N$ such that $\ds\liminf_{n \ra \infty}\al_n(1-\al_n) > 0$. Then, the sequence $\{x_n\}_{n\geq 0}$ converges to some element of $F$.
\end{cor}

\noindent
If, in Theorem 2, we have that $f_1 \equiv f_2 \equiv 0,$ then we have the following corollary.

\begin{cor}
Let $C$ be a nonempty closed convex subset of a uniformly convex and uniformly smooth real Banach space $E$. Let $B_1, B_2 : C \ra E^*$ be continuous monotone mappings. Let $T_j : C \ra C, \; j = 1, 2, \ldots, d$ be a finite family of continuous relatively nonexpansive mappings. Let $F := \ds\left[\bigcap_{j=1}^dF(T_j)\right]\cap\left[\bigcap_{k=1}^2VI(B_k, C)\right] \neq \emptyset$. Let $\{x_n\}_{n\geq 0}$ be a sequence defined by

\begin{equation}
\left\{ \begin{array}{ll}
x_0 \in C_0 = C, \quad \mbox{chosen arbitrarily},\\
y_n = J^{-1}(\alpha_{n}Jx_n+(1-\al_n)JT_{n+1}x_n); u_n, v_n \in C \; \mbox{s.t.}\\
\ltl B_1y_n, y-u_n\rtl+\frac{1}{r_n}\ltl y-u_n, Ju_n-Jy_n\rtl \geq 0
, \forall \; y \in C,\\
\ltl B_2y_n, y-v_n\rtl+\frac{1}{r_n}\ltl y-v_n, Jv_n-Jy_n\rtl \geq 0
, \forall \; y \in C,\\
w_n = J^{-1}(\beta Ju_n+(1-\beta) Jv_n),\\
C_{n+1} = \{z \in C_n: \; \phi(z, w_n) \leq \phi(z, x_n)\},\\
x_{n+1} = \Pi_{C_{n+1}}(x_0), \; n\geq 0,
\end{array}
\right.
\end{equation}

\noindent
where $T_n = T_{n(\bmod m)}$ and $J$ is the normalized duality mapping on $E$; $\{r_n\}_{n \geq 1} \su [c_1, \infty)$ for some $c_1 > 0$, $\beta$, $\al_n \in (0,1)$ for all $n \in \N$ such that $\ds\liminf_{n \ra \infty}\al_n(1-\al_n) > 0$. Then, the sequence $\{x_n\}_{n\geq 0}$ converges to some element of $F$.
\end{cor}

\noindent
If $E = H$, a Hilbert space, then $E$ is $2$-uniformly convex and uniformly smooth real Banach space, in this case, continuous relatively nonexpansive mapping reduces to continuous quasi-nonexpansive mapping. Furthermore, $J = I$, identity operator on $H$ and $\Pi_C = P_C$, projection mapping from $H$ into $C$. Thus, the following corollaries hold.

\begin{cor}
Let $C$ be a nonempty closed convex subset of a real Hilbert space $H$. Let $f_1, f_2 : C\times C \ra \R, \; k = 1,2, \cdots, q$ be bifunctions satisfying (A1) - (A4) and $B_1, B_2 : C \ra E^*$ be monotone continuous mappings. Let $T_j : C \ra C, \; j = 1, 2, \cdots,d$ be finite family of continuous quasi-nonexpansive mappings and $A_i : C \ra E^*, \; i = 1,2, \cdots, m$ be finite family of $\ga_i$-inverse strongly monotone operators with constants $\ga_i \in (0,1), \; i = 1, 2, \cdots, m$. Let $F := \ds\left[\bigcap_{j=1}^dF(T_j)\right]\cap\left[\bigcap_{i=1}^mA_i^{-1}(0)\right]\cap\left[\bigcap_{k=1}^2GEP(f_k, B_k)\right] \neq \emptyset$. Let $\{x_n\}_{n\geq 0}$ be a sequence defined by

\begin{equation}
\left\{ \begin{array}{ll}
x_0 \in C_0 = C, \quad \mbox{chosen arbitrarily},\\
z_n = P_{C}(x_n - \la_nA_{n+1}x_n)),\\
y_n = \alpha_{n}x_n+(1-\al_n)T_{n+1}z_n; u_n, v_n \in C \; \mbox{s.t.}\\
f_1(u_n,y)+\ltl B_1y_n, y-u_n\rtl+\frac{1}{r_n}\ltl y-u_n, u_n-y_n\rtl \geq 0
, \forall \; y \in C,\\
f_2(v_n,y)+\ltl B_2y_n, y-v_n\rtl+\frac{1}{r_n}\ltl y-v_n, v_n-y_n\rtl \geq 0
, \forall \; y \in C,\\
w_n = \beta u_n+(1-\beta) v_n,\\
C_{n+1} = \{z \in C_n: \; \lt z-w_n\rt \leq \lt z-x_n\rt\},\\
x_{n+1} = P_{C_{n+1}}(x_0), \; n\geq 0,
\end{array}
\right.
\end{equation}

\noindent
where $A_n = A_{n(\bmod m)}$, $T_n = T_{n(\bmod m)}$ and $\{r_n\}_{n \geq 1} \su [c_1, \infty)$ for some $c_1 > 0$, $\beta$, $\al_n \in (0,1)$ for all $n \in \N$ such that $\ds\liminf_{n \ra \infty}\al_n(1-\al_n) > 0$; and $\{\la_n\}_{n\geq 1}$ is a sequence in $[a, b]$ for some $0<a<b<\frac{c^2\ga}{2}$, where $c$ is the $2$-uniformly convex constant of $E$ and $\ga = \ds\min_{1\leq i\leq m}\ga_i$. Then, the sequence $\{x_n\}$ converges strongly to a point of $F$.
\end{cor}

\begin{cor}
Let $C$ be a nonempty closed convex subset of a real Hilbert space $H$. Let $f_1, f_2 : C\times C \ra \R, \; k = 1,2, \cdots, q$ be bifunctions satisfying (A1) - (A4) and $B_1, B_2 : C \ra E^*$ be monotone continuous mappings. Let $T_j : C \ra C, \; j = 1, 2, \cdots,d$ be finite family of continuous quasi-nonexpansive mappings. Let $F := \ds\left[\bigcap_{j=1}^dF(T_j)\right]\cap\left[\bigcap_{k=1}^2GEP(f_k, B_k)\right] \neq \emptyset$. Let $\{x_n\}_{n\geq 0}$ be a sequence defined by

\begin{equation}
\left\{ \begin{array}{ll}
x_0 \in C_0 = C, \quad \mbox{chosen arbitrarily},\\
y_n = \alpha_{n}x_n+(1-\al_n)T_{n+1}x_n; u_n, v_n \in C \; \mbox{s.t.}\\
f_1(u_n,y)+\ltl B_1y_n, y-u_n\rtl+\frac{1}{r_n}\ltl y-u_n, u_n-y_n\rtl \geq 0
, \forall \; y \in C,\\
f_2(v_n,y)+\ltl B_2y_n, y-v_n\rtl+\frac{1}{r_n}\ltl y-v_n, v_n-y_n\rtl \geq 0
, \forall \; y \in C,\\
w_n = \beta u_n+(1-\beta) v_n,\\
C_{n+1} = \{z \in C_n: \; \lt z-w_n\rt \leq \lt z-x_n\rt\},\\
x_{n+1} = P_{C_{n+1}}(x_0), \; n\geq 0,
\end{array}
\right.
\end{equation}

\noindent
where $T_n = T_{n(\bmod m)}$ and $\{r_n\}_{n \geq 1} \su [c_1, \infty)$ for some $c_1 > 0$, $\beta$, $\al_n \in (0,1)$ for all $n \in \N$ such that $\ds\liminf_{n \ra \infty}\al_n(1-\al_n) > 0$. Then, the sequence $\{x_n\}$ converges strongly to a point of $F$.
\end{cor}

\begin{cor}
Let $C$ be a nonempty closed convex subset of a real Hilbert space $H$. Let $T_j : C \ra C, \; j = 1, 2, \cdots,d$ be finite family of continuous quasi-nonexpansive mappings. Let $B_1, B_2 : C \ra E^*$ be monotone continuous mappings. Let $F := \ds\left[\bigcap_{j=1}^dF(T_j)\right]\cap\left[\bigcap_{k=1}^2VI(B_k, C)\right] \neq \emptyset$. Let $\{x_n\}_{n\geq 0}$ be a sequence defined by

\begin{equation}
\left\{ \begin{array}{ll}
x_0 \in C_0 = C, \quad \mbox{chosen arbitrarily},\\
y_n = \alpha_{n}x_n+(1-\al_n)T_{n+1}x_n; u_n, v_n \in C \; \mbox{s.t.}\\
\ltl B_1y_n, y-u_n\rtl+\frac{1}{r_n}\ltl y-u_n, u_n-y_n\rtl \geq 0
, \forall \; y \in C,\\
\ltl B_2y_n, y-v_n\rtl+\frac{1}{r_n}\ltl y-v_n, v_n-y_n\rtl \geq 0
, \forall \; y \in C,\\
w_n = \beta u_n+(1-\beta) v_n,\\
C_{n+1} = \{z \in C_n: \; \lt z-w_n\rt \leq \lt z-x_n\rt\},\\
x_{n+1} = P_{C_{n+1}}(x_0), \; n\geq 0,
\end{array}
\right.
\end{equation}

\noindent
where $T_n = T_{n(\bmod m)}$ and $\{r_n\}_{n \geq 1} \su [c_1, \infty)$ for some $c_1 > 0$, $\beta$, $\al_n \in (0,1)$ for all $n \in \N$ such that $\ds\liminf_{n \ra \infty}\al_n(1-\al_n) > 0$. Then, the sequence $\{x_n\}$ converges strongly to a point of $F$.
\end{cor}

\begin{cor}
Let $C$ be a nonempty closed convex subset of a real Hilbert space $H$. Let $f_1, f_2 : C\times C \ra \R, \; k = 1,2, \cdots, q$ be bifunctions satisfying (A1) - (A4). Let $T_j : C \ra C, \; j = 1, 2, \cdots,d$ be finite family of continuous quasi-nonexpansive mappings. Let $F := \ds\left[\bigcap_{j=1}^dF(T_j)\right]\cap\left[\bigcap_{k=1}^2EP(f_k)\right] \neq \emptyset$. Let $\{x_n\}_{n\geq 0}$ be a sequence defined by

\begin{equation}
\left\{ \begin{array}{ll}
x_0 \in C_0 = C, \quad \mbox{chosen arbitrarily},\\
y_n = \alpha_{n}x_n+(1-\al_n)T_{n+1}x_n; u_n, v_n \in C \; \mbox{s.t.}\\
f_1(u_n,y)+\frac{1}{r_n}\ltl y-u_n, u_n-y_n\rtl \geq 0
, \forall \; y \in C,\\
f_2(v_n,y)+\frac{1}{r_n}\ltl y-v_n, v_n-y_n\rtl \geq 0
, \forall \; y \in C,\\
w_n = \beta u_n+(1-\beta) v_n,\\
C_{n+1} = \{z \in C_n: \; \lt z-w_n\rt \leq \lt z-x_n\rt\},\\
x_{n+1} = P_{C_{n+1}}(x_0), \; n\geq 0,
\end{array}
\right.
\end{equation}

\noindent
where $T_n = T_{n(\bmod m)}$ and $\{r_n\}_{n \geq 1} \su [c_1, \infty)$ for some $c_1 > 0$, $\beta$, $\al_n \in (0,1)$ for all $n \in \N$ such that $\ds\liminf_{n \ra \infty}\al_n(1-\al_n) > 0$. Then, the sequence $\{x_n\}$ converges strongly to a point of $F$.
\end{cor}

\noindent
Note that, if, in Theorem 1, we replace\\
$``F := \ds\left[\bigcap_{j=1}^dF(T_j)\right]\cap\left[\bigcap_{i=1}^mA_i^{-1}(0)\right]\cap\left[\bigcap_{k=1}^2GEP(f_k, B_k)\right] \neq \emptyset"$ with $``\lt Ax\rt \leq \lt Ax-Ap\rt \; \forall \; p \in F,\;\\ x \in C \; \mbox{and} \; F := \ds\left[\bigcap_{j=1}^dF(T_j)\right]\cap\left[\bigcap_{i=1}^mVI(A_i, C)\right]\cap\left[\bigcap_{k=1}^2GEP(f_k, B_k)\right] \neq \emptyset"$,  then we get the following theorem.

\begin{theo}
Let $C$ be a nonempty closed convex subset of $2$-uniformly convex and uniformly smooth real Banach space $E$. Let $f_1, f_2 : C \times C \ra \R$ be bifunctions satisfying (A1) - (A4) and $B_1, B_2 : C \ra E^*$ be continuous monotone mappings. Let $T_j : C \ra C, \; j = 1, 2, \ldots, d$ be a finite family of continuous relatively nonexpansive mappings and $A_i : C \ra E^*, \; i = 1, 2, \ldots, m$ be a finite family of $\ga_i$-inverse strongly monotone operators with constants $\ga_i \in (0,1), \; i = 1, 2, \ldots, m.$ Let $F := \ds\left[\bigcap_{j=1}^dF(T_j)\right]\cap\left[\bigcap_{i=1}^mVI(A_i, C)\right]\cap\left[\bigcap_{k=1}^2GEP(f_k, B_k)\right] \neq \emptyset$. Let $\{x_n\}_{n\geq 0}$ be a sequence defined by (13). Then, the sequence $\{x_n\}_{n\geq 0}$ converges to some element of $F$.
\end{theo}
\begin{prf}
Let $p \in F$. Then by assumption $\lt Ax\rt \leq \lt Ax - Ap\rt \; \forall \; x \in C$ and in particular \\ $\lt Ax\rt \leq \lt Ax - Ap\rt = 0$ which implies that $Ap = 0$ and so $p \in A^{-1}(0)$. Therefore, the conclusion follows from Theorem 1.
\end{prf}

\noindent
The proof of the following theorem can be easily obtained from the method of proof of Theorem 1.

\begin{theo}
Let $C$ be a nonempty closed convex subset of $2$-uniformly convex and uniformly smooth real Banach space $E$. Let $f_k : C \times C \ra \R \; k = 1,2,\ldots, q$ be finite family of bifunctions satisfying (A1) - (A4) and $B_k : C \ra E^* \; k = 1,2,\ldots, q$ be finite continuous monotone mappings. Let $T_j : C \ra C, \; j = 1, 2, \ldots, d$ be a finite family of continuous relatively nonexpansive mappings and $A_i : C \ra E^*, \; i = 1, 2, \ldots, m$ be a finite family of $\ga_i$-inverse strongly monotone operators with constants $\ga_i \in (0,1), \; i = 1, 2, \ldots, m.$ Let $F := \ds\left[\bigcap_{j=1}^dF(T_j)\right]\cap\left[\bigcap_{i=1}^mVI(A_i, C)\right]\cap\left[\bigcap_{k=1}^qGEP(f_k, B_k)\right] \neq \emptyset$. Let $\{x_n\}_{n\geq 0}$ be a sequence defined by

\begin{equation}
\left\{ \begin{array}{ll}
x_0 \in C_0 = C, \quad \mbox{chosen arbitrarily},\\
z_n = \Pi_{C}J^{-1}(Jx_n - \la_nA_{n+1}x_n)),\\
y_n = J^{-1}(\alpha_{n}Jx_n+(1-\al_n)JT_{n+1}z_n), \; u_{k,n} \in C, \; k = 1, 2, \ldots, q \; \mbox{s.t.}\\
f_1(u_{1,n},y)+\ltl B_1y_n, y-u_{1,n}\rtl+\frac{1}{r_n}\ltl y-u_{1,n}, Ju_{1,n}-Jy_n\rtl \geq 0
, \forall \; y \in C,\\
f_2(u_{2,n},y)+\ltl B_2y_n, y-u_{2,n}\rtl+\frac{1}{r_n}\ltl y-u_{2,n}, Ju_{2,n}-Jy_n\rtl \geq 0
, \forall \; y \in C,\\
\vdots \\
f_q(u_{q,n},y)+\ltl B_qy_n, y-u_{q,n}\rtl+\frac{1}{r_n}\ltl y-u_{q,n}, Ju_{q,n}-Jy_n\rtl \geq 0
, \forall \; y \in C,\\
w_n = J^{-1}(\beta_1 Ju_{1,n}+\beta_2 Ju_{2,n}+\ldots +\beta_q Ju_{q,n}),\\
C_{n+1} = \{z \in C_n: \; \phi(z, w_n) \leq \phi(z, x_n)\},\\
x_{n+1} = \Pi_{C_{n+1}}(x_0), \; n\geq 0,
\end{array}
\right.
\end{equation}

\noindent
where $A_n = A_{n(\bmod m)}$, $T_n = T_{n(\bmod m)}$ and J is the normalized duality mapping on $E$; $\{r_n\}_{n \geq 1} \su [c_1, \infty)$ for some $c_1 > 0$, $\beta$, $\al_n \in (0,1)$ for all $n \in \N, \; k = 1,2, \dots, q$ such that $\ds\sum_{k=1}^{q}\beta_k, \; \ds\liminf_{n \ra \infty}\al_n(1-\al_n) > 0$; and $\{\la_n\}_{n\geq 1}$ is a sequence in $[a, b]$ for some $0<a<b<\frac{c^2\ga}{2}$, where $c$ is the $2$-uniformly convex constant of $E$ and $\ga = \ds\min_{1\leq i\leq m}\ga_i$. Then, the sequence $\{x_n\}$ converges strongly to a point of $F$.
\end{theo}

\subsection*{\normalsize \center CONCLUSION}
\noindent
Corollary 1 improves Theorem 3.1 of Takahashi and Zembayashi \cite{takahashi2008strong} to a finite family of relatively nonexpansive mappings and equilibrium problems. Theorem 3.1 of Li and Su \cite{hongyu2010strong} is a special case of Theorem 3 in which $i = 1, \; k=1, \; B_1 \equiv B_2 \equiv 0$ and $T_j \equiv I$ for all $j = 1, \ldots, d$. Slight modifications of the iteration schemes studied in this paper extend the results of \textit{Ofoedu et al}. \cite{ofoedu6} from Hilbert space to $2$-uniformly convex Banach space. Our theorems improve and generalize the main results of Takahashi and Zembayashi \cite{takahashi2008strong}, Li and Su \cite{hongyu2010strong}, \textit{Ofoedu et al}. \cite{ofoedu6} and several other results which are announced recently. Our iteration process, method of proof and corollaries are of independent interest.

\subsection*{\normalsize DECLARATION AND ACKNOWLEDGMENT}
\noindent
We wish to declare that all authors contributed equally to the preparation of this paper. We thank the unanimous referee(s) for critical review of the manuscript and for giving suggestions that helped improve the paper.

\renewcommand{\bibname}{References}
\bibliographystyle{plain}
\bibliography{references}

\end{document}